\documentclass[12pt]{amsart}

\usepackage{xy}
\usepackage{amscd,amssymb,latexsym}
\usepackage{graphicx}
\usepackage{amsfonts}

\input{xy}
\xyoption{all}

\newtheorem*{corollary*}{Corollary}

\newtheorem*{proposition*}{Proposition}
\newtheorem*{propositionA}{Proposition A}
\newtheorem{lemma}{Lemma}
\newtheorem*{theorem*}{Theorem}

\theoremstyle{definition}
\newtheorem*{definition*}{Definition}

\newtheorem*{remark*}{Remark}

\newcommand{\comment}[1]{}

\title[Surfaces Violating BMY in Positive Characteristic]{Surfaces Violating Bogomolov-Miyaoka-Yau in Positive Characteristic}
\author{Robert W. Easton}
\begin{document}

\begin{abstract}
The Bogomolov-Miyaoka-Yau inequality asserts that the Chern numbers
of a surface $X$ of general type in characteristic 0 satisfy the
inequality $c_1^2\leq 3c_2$, a consequence of which is
$\frac{K_X^2}{\chi ({\mathcal O}_X)}\leq 9$.  This inequality  fails in
characteristic $p$, and here we produce infinite families of
counterexamples for large $p$.  Our method parallels a construction
of Hirzebruch, and relies on a construction of abelian covers due to
Catanese and Pardini.
\end{abstract}

\maketitle

\section{Introduction}

In 1976, Bogomolov \cite{B} proved that if $X$ is a complex surface
of general type, then the Chern numbers of $X$ satisfy the
inequality $c_1^2\leq 4c_2$, where $c_i=c_i(\mathcal{T}_X)$.  This
inequality was improved independently by Miyaoka \cite{M2} and Yau
\cite{Y} to $c_1^2\leq 3c_2$, which is sharp by examples of
Hirzebruch \cite{HA}--\cite{HC}. A consequence of this latter
inequality is that if $X$ is a surface of general type in
characteristic $0$, then $\frac{K_X^2}{\chi ({\mathcal O}_X)}\leq 9$. This
inequality does not hold in positive characteristic, as can already
be seen in families of surfaces constructed by Lang \cite{L}. In this paper
we produce infinite families of counterexamples in large
characteristics.

The surfaces we generate parallel a construction of Hirzebruch.  In
\cite{HB}, Hirzebruch produced smooth complex surfaces of general
type with $c_1^2=3c_2$ by constructing abelian covers of
$\mathbb{P}^2$ branched over specific configurations of lines.  Here
we apply a similar procedure in positive characteristic, producing
surfaces branched over pathological configurations of lines unique
to positive characteristic.

More explicitly, following the notation of Vakil \cite[\S4]{V}, fix a prime
$q$ and positive integer $n$, and let
$G=(\mathbb{Z}/q\mathbb{Z})^n$. Let $\langle \cdot,\cdot \rangle: G
\times G^{\vee} \to \mathbb{Z}/q\mathbb{Z}$ be the natural pairing,
which we extend in the obvious way to a map $\langle \cdot,\cdot
\rangle: G \times G^{\vee} \to \mathbb{Z}$. (Note that this map is no
longer bilinear, but that shall not be needed here.) Suppose $k$ is a
field of characteristic $p\neq q$, and let $S$ be the blowup of
$\mathbb{P}_k^2$ at the $(p^2+p+1)$ $\mathbb{F}_p$-valued points.

\bigskip
\pagebreak
\noindent Suppose we have two maps:

\smallskip
$D: G \to {\rm Div}(S)$

$L: G^{\vee} \to {\rm Pic}(S)$
\smallskip

\noindent We say $(D,L)$ satisfy the \emph{cover condition} if $D_0=0$ and if, for all $\chi \in G^{\vee}$,
\[
qL_\chi=\sum_{\sigma \in G} \langle \chi, \sigma \rangle D_\sigma,
\]
with equality taking place in ${\rm Pic}(S)$.

We then have the following result, due to Catanese \cite{C} and Manetti \cite[\S3]{M1} in the case $G=(\mathbb{Z}/2\mathbb{Z})^n$ and Pardini \cite[\S4]{P} in general:

\begin{propositionA}
Suppose $(D,L)$ satisfy the cover condition, and
\begin{enumerate}

\item[i)] $D_\sigma$ are all nonsingular curves;

\item[ii)] no three meet at a point; and

\item[iii)] if $D_\sigma, D_{\sigma'}$ meet, then they are transverse, and
$\sigma, \sigma'$ are independent in $G$.
\end{enumerate}

Then:
\begin{enumerate}

\item[(1)] There exists a $G$-cover $\pi: \tilde S \to S$ with \~{S} a
nonsingular surface, and with branch locus $D = \cup D_{\sigma}$;

\item[(2)] $q^n K_{\tilde S}=\pi^*(q^n K_S+q^{n-1}(q-1)\sum_{\sigma}D_{\sigma})$, with equality holding in ${\rm Pic}(\tilde S)$;

\item[(3)] $K_{\tilde S}^2=q^n\left( K_S+\frac{q-1}{q}\sum_{\sigma}D_{\sigma}\right)^2$; and

\item[(4)] $\chi ({\mathcal O}_{\tilde S})=q^n\chi ({\mathcal O}_S)+\frac{1}{2}\sum_{\chi}L_{\chi}\cdot (L_{\chi}+K_S)$.
\end{enumerate}
\end{propositionA}

Our general strategy is to build surfaces $\tilde S$ using divisors
on $S$ which exploit pathologies in positive characteristic. One
such pathology is the following:  If we consider the set of lines
through pairs of $\mathbb{F}_p$-valued points in $\mathbb{P}_k^2$,
we find there are more such lines than ``would be expected" in the
characteristic $0$ case.  To exploit this pathology, we let $C$ be
the union of these lines and take the strict transform $\tilde{C}$
of $C$ to be one of our divisors.  We then choose the remaining
divisors to be unions of appropriate numbers of pullbacks of general lines in $\mathbb{P}^2$.  With particularly cogent choices of divisors, line bundles
$L_{\chi}$ can be defined directly from the covering condition.

Our main result is the following:

\begin{theorem*}For each prime $q\geq 3$ and all sufficiently large primes $p\equiv -1 \pmod q$, there exists a surface $\tilde S$ in characteristic $p$, nonsingular and of general type, with $\frac{K_{\tilde S}^2}{\chi ({\mathcal O}_{\tilde S})}>9$.
\end{theorem*}

\paragraph{\bf{Acknowledgements:}}

I thank R. Vakil for suggesting this strategy of constructing
interesting surfaces in positive characteristic, as well as
providing insight throughout the process. His guidance has been
invaluable.

I also thank R. Pardini for her many helpful suggestions,  and especially for discovering a serious error in the original draft of this paper.

\section{The Construction}
Let $p,q$ be primes with $q\geq 3$ and $p\equiv-1 \pmod q$, and let $n\geq 3$ be any integer.

For each $\sigma \in G=(\mathbb{Z}/q\mathbb{Z})^n$, write $\sigma = (\sigma_1,\ldots,\sigma_n),$ with $\sigma_i \in \{ 0,\ldots,q-1 \}$.  Let $\vec{e}_1,\ldots,\vec{e}_n$ represent the standard unit `vectors'. Also, for each $\gamma \in G$, let $\chi_{\gamma} \in G^{\vee}$ be defined by $\langle \chi_{\gamma},\vec{e}_i
\rangle=\gamma_i$.  Throughout, let $k$ be a field with ${\rm char}(k)=p$.

Observe that there are $(p^2+p+1)$ $\mathbb{F}_p$-valued points in
$\mathbb{P}_k^2,$ and the same number of such lines. Also observe
that through any fixed $\mathbb{F}_p$-valued point there pass
exactly $p+1$ of these lines, and each line contains exactly $p+1$
of the $\mathbb{F}_p$-valued points.  Let $S$ be the blowup of
$\mathbb{P}_k^2$ at the $\mathbb{F}_p$-valued points, and let
$\tilde{C}$ be the strict transform of the union $C$ of the lines
through these points.  Note that $K_S\sim
-3H+\sum_{i=1}^{p^2+p+1}E_i$, and $[\tilde{C}]\sim
(p^2+p+1)H-(p+1)\sum E_i$, where $H$ is the hyperplane class of $S$
and the $E_i$ are the exceptional divisors.

We will set $D_{\vec{e}_1}=[\tilde{C}]$, and choose the remaining divisors $D_{\sigma}$ to be the union of the pullback of general lines in $\mathbb{P}^2$.  To ensure that the covering condition is `solvable', some care is needed in choosing the number of sections for each divisor. To that end, let $F\subset G$ be any subset such that $\{\vec{0},\vec{e}_1+\vec{e}_2,\vec{e}_1,\ldots,\vec{e}_n\} \subset F,$ $F\cap kF=\{\vec{0}\}$ for $k=2,\ldots,q-1,$ and $G=\bigcup_{k=1}^{q-1} kF$.
More concretely, $F$ should contain the required elements and exactly one representative from each distinct subset of the form $\{k\sigma : k=1,\ldots,q-1\}.$

\begin{lemma} There exists a collection of integers $m_{\sigma}\in \{0,\ldots,q-1\}$ such that
\begin{enumerate}
\item[(i)]$m_{\vec{e}_1}=1$;
\item[(ii)]$m_{\sigma}=0$, for all $ \sigma \notin F$; and
\item[(iii)]$\sum_{\sigma \in G}m_{\sigma}\sigma=\vec{0}$.
\end{enumerate}
\end{lemma}

\begin{proof}
Define the set $\tilde F=F\backslash\{\vec{0},\vec{e}_1+\vec{e}_2,\vec{e}_1,\ldots,\vec{e}_n\}$.
Let $m_{\sigma}=0$ for all $\sigma \notin F$, $m_{\sigma}=1$ for all $\sigma \in \tilde F$, $m_{\vec{0}}=0$, and $m_{\vec{e}_1}=1$.
Consider the element
\[
(A_1,\ldots,A_n):=\sum_{\sigma \in \tilde F}\sigma \in G,
\]
Working modulo $q$, choose $m_{\vec{e}_1+\vec{e}_2}\equiv -A_1-1$, $m_{\vec{e}_2}\equiv -A_2+A_1+1$, and $m_{\vec{e}_i}\equiv -A_i$ for each $i=3,\ldots,n$.  The result follows.
\end{proof}

Let $\mathfrak{S}$ be any collection of such integers.  We're then in a position to choose our divisors and line bundles.
For a fixed choice of $r_{\sigma}\in \mathbb{Z}_{\geq 0}$, let
\begin{itemize}
\item $D_{\vec{0}}=0$
\item $D_{\vec{e_1}}=[\tilde{C}]$
\item $D_{\sigma}=\left\{
                \begin{array}{ll}
                    \text{union of }qr_{\sigma}+m_{\sigma}\text{ general sections of }S,& \text{if} \; \sigma \in F\backslash\{\vec{0},\vec{e_1}\}\\
                    0,& \text{otherwise.}
                \end{array}
                \right.$
\end{itemize}
Note that
\begin{itemize}
\item $D_{\vec{e_1}}\sim(p^2+p+1)H-(p+1)\sum_{i=1}^{p^2+p+1}E_i$
\item $D_{\sigma}\sim\left\{
                \begin{array}{ll}
                  (qr_{\sigma}+m_{\sigma})H,& \text{if} \; \sigma \in F\backslash\{\vec{0},\vec{e_1}\}\\
                   0,& \text{otherwise.}
                \end{array}
                \right.$
\end{itemize}

Our choice of the sets  $F, \mathfrak{S}$ ensures that conditions
(i)-(iii) of Proposition A are satisfied.  Moreover, from the lemma
(and the fact that $p\equiv -1 \pmod q$) it is straightforward
to verify that $q|\sum_{\sigma \in G} \langle \chi_{\gamma},\sigma
\rangle D_{\sigma}$, for every $\gamma \in G.$  We may therefore
define
\[
L_{\chi_{\gamma}}:=\frac{1}{q}\sum_{\sigma \in G}\langle
\chi_{\gamma},\sigma \rangle D_{\sigma},
\]
for each $\gamma \in G,$ which clearly satisfy the covering
condition.  For later reference, these choices of $D_{\sigma}$ and
$L_{\chi}$ give

\begin{equation}
\sum_{\sigma \in G}D_{\sigma} \sim \left( p^2+p+1+\sum_{\sigma \in
F\backslash\{\vec{0},\vec{e_1}\}}(qr_{\sigma}+m_{\sigma})\right)H-(p+1)\sum_{i=1}^{p^2+p+1}E_i
\label{eq:1}
\end{equation}

\noindent{and}

\begin{equation}
\begin{split}
\sum_{\gamma \in G}L_{\chi_{\gamma}}&\cdot \left(L_{\chi_{\gamma}}+K_S\right)=\\
\frac{1}{q^2}\sum_{\gamma \in G}[&(\gamma_1(p^2+p+1)+f(\gamma,\sigma))\cdot(-3q+\gamma_1(p^2+p+1)+f(\gamma,\sigma))\\
&+\gamma_1(p+1)(q-\gamma_1(p+1))(p^2+p+1)],
\end{split}
\label{eq:2}
\end{equation}

\noindent where $f(\gamma,\sigma)=\sum_{\sigma \in F\backslash\{\vec{0},\vec{e_1}\}}\overline{\gamma \cdot \sigma}(qr_{\sigma}+m_{\sigma})$,
and where $\overline{\gamma \cdot \sigma}$ is the reduced residue of $\gamma \cdot \sigma$ modulo $q$.

Before continuing, we should check that the surface $\tilde S$ constructed from these $(D,L)$ is of general type.

\begin{lemma}
Let $\tilde S$ be the abelian cover of $S$ with building data $(D,L)$ given above.
Then $\tilde S$ is of general type.
\end{lemma}

\begin{proof}
Recall that $\tilde{S}$ is of general type if and only if the
divisor $K_{\tilde{S}}$ is big (and that on a smooth projective
surface, a divisor $D$ is big if for some integer $n>0$ one has
$nD\sim E+A$ for some effective divisor $E$ and ample divisor $A$).
Since $q^n K_{\tilde S}=\pi^{\ast}(q^n
K_S+q^{n-1}(q-1)\sum_{\sigma}D_{\sigma})$ and $\pi$ is a finite
morphism, it suffices to prove that the divisor $q^n
K_S+q^{n-1}(q-1)\sum_{\sigma}D_{\sigma}$ is big. By our choice of
divisors $D_{\sigma}$, we have
\begin{align*}
q^n & K_S+q^{n-1}(q-1) \sum_{\sigma}D_{\sigma}\sim q^n\left(-3H + \sum_i E_i\right)\\
&\quad +q^{n-1}(q-1)\left( \left( p^2+p+1+\sum_{\sigma \in F\backslash \{\vec{0},\vec{e_1}\}} (qr_{\sigma}+m_{\sigma}) \right)H - (p+1)\sum_i E_i \right)\\
&\sim q^{n-1}\left(-3q+(q-1)\left( p^2+p+1+\sum_{\sigma \in F\backslash \{\vec{0},\vec{e_1}\}} (qr_{\sigma}+m_{\sigma})\right) \right)H\\
&\quad+q^{n-1}\left( q-(q-1)(p+1)\right)\sum_i E_i\\
&\sim q^{n-1}(q-1)\left((p^2+p+1)H-(p+1)\sum_i E_i\right)\\
&\quad+q^{n-1}\left( (q-1)\sum_{\sigma \in F\backslash
\{\vec{0},\vec{e_1}\}} (qr_{\sigma}+m_{\sigma})-3q\right)H+q^n\sum_i
E_i
\end{align*}

\noindent Now, since $r_{\sigma},m_{\sigma}\geq 0$ for all $\sigma$
and $m_{\sigma}=1$ for all $\sigma \in \tilde F$, we have
\[
\sum_{\sigma \in
F\backslash\{\vec{0},\vec{e}_1\}}(qr_{\sigma}+m_{\sigma})\geq
|\tilde F| =\frac{q^n-1}{q-1}-n,
\]
and so
\[
(q-1)\sum_{\sigma \in F\backslash \{\vec{0},\vec{e_1}\}}
(qr_{\sigma}+m_{\sigma})-3q \geq q^n -1-(q-1)n-3q=q^n-(n+3)q+n-1 >0,
\]
(for $q,n \geq 3$).  We can therefore write
\[
q^n K_S+q^{n-1}(q-1) \sum_{\sigma}D_{\sigma}\sim q^{n-1}(q-1)[\tilde
C]+mH+q^n\sum E_i
\]
for some $m > 0$.  Since $q^{n-1}(q-1)[\tilde C] +q^n\sum E_i$ is
effective and $H$ is ample, the result follows.
\end{proof}

\section{Results}

Using equations \eqref{eq:1} and \eqref{eq:2} above, it is
straightforward to calculate $\frac{K^2}{\chi}$ for given values of
$q,n$ and $p$. For simplicity, for the remainder of this paper we'll
assume $r_{\sigma}=1$ for all $\sigma \in G.$  With this convention
and the above notation, let us define
\[
R(q,n,p;F,\mathfrak{S}):=\frac{K_{\tilde S}^2}{\chi ({\mathcal O}_{\tilde S})}.
\]
This ratio will clearly depend on the particular choice of the sets
$F$ and $\mathfrak{S}$. As we will see, however, our main result
will be independent of these choices.

Before proceeding to the main result, let us first give a few
explicit calculations, as examples of the types of ratios one
obtains. In the following examples, the set $F$ was canonically
chosen (via a construction based on the natural ordering of elements
in $(\mathbb{Z}/q\mathbb{Z})^n$), and the set $\mathfrak{S}$ was
then consistently chosen as in the proof of Lemma 1.  (It is not worth going into the specifics of the construction, since it is used only in the following examples and will not be needed for any of the later results.)

For $n=3$ we found:

\begin{align*}
R(3,3,p)&=\frac{3(7568+351p+351p^2+8p^3)}{3023+134p+134p^2+2p^3}, \\
R(5,3,p)&=\frac{497024+5647p+5647p^2+24p^3}{63197+708p+708p^2+2p^3}, \\
R(7,3,p)&=\frac{7112888+32015p+32015p^2+48p^3}{896249+4006p+4006p^2+4p^3}. \\
\end{align*}

As another example, for $n=4$ we found:

\begin{align*}
R(3,4,p)&=\frac{3(91808+1215p+1215p^2+8p^3)}{35045+458p+458p^2+2p^3}, \\
R(5,4,p)&=\frac{13727024+29647p+29647p^2+24p^3}{1721447+3708p+3708p^2+2p^3}, \\
R(7,4,p)&=\frac{365995160+229583p+229583p^2+48p^3}{45800437+28702p+28702p^2+4p^3}.
\end{align*}

More relevant than the ratios themselves, however, will be their behavior as we let $p$ tend to infinity.
So, let us define
\[
R(q,n):= \lim_{p \rightarrow \infty}R(q,n,p).
\]

\noindent From our above list, we then see that

\medskip
$R(3,3)= R(3,4)= R(5,3)= R(5,4)=R(7,3)= R(7,4)=12$

\medskip
Let us now prove the following:

\begin{proposition*}For all primes $q \geq 3,$ and integers $n\geq 3,$ \[R(q,n)=12\]
independent of the choices of the sets $F, \mathfrak{S}$.
\end{proposition*}

The bulk of the proof is contained in the following lemma.

\begin{lemma} For fixed $q$ and $n$, we have the following estimates, independent of the choices of the sets $F, \mathfrak{S}$:
\begin{enumerate}
\item[i)] $K_S^2=O(p^2)$;
\item[ii)] $K_S \cdot \sum_{\sigma}D_{\sigma} = p^3+O(p^2)$;
\item[iii)] $\left( \sum_{\sigma} D_{\sigma} \right)^2 = -p^3+O(p^2)$;
\item[iv)] $\sum_{\chi} L_{\chi} \cdot (L_{\chi}+K_S) = \left(\frac{1}{6}q^{n-2}(q^2-1)\right)p^3+O(p^2)$.
\end{enumerate}
\end{lemma}

\begin{proof}  First note that the sizes of the sets $F$ and $\mathfrak{S}$ are bounded in terms of $q$ and $n$, as are also the elements of $\mathfrak{S}$.
Thus, for fixed $q$ and $n$, sums over $F$ of elements in $\mathfrak{S}$ contribute $O(1)$.  This observation considerably simplifies many of the following estimates.

\smallskip
\noindent (i)  Since $K_S=-3H+\sum_{i=1}^{p^2+p+1}E_i$, we
immediately have $K_S^2=9-(p^2+p+1)=O(p^2)$.  For the remaining
estimates, we rely heavily on equations \eqref{eq:1} and
\eqref{eq:2}.

\smallskip
\noindent (ii) Observe that
\begin{align*}
K_S \cdot \sum_{\sigma \in G} D_{\sigma} &= K_S \cdot \left( \left( p^2+p+1+\sum_{\sigma \in F\backslash\{\vec{0},\vec{e}_1\}}(q+m_{\sigma})\right)H-(p+1)\sum_{i=1}^{p^2+p+1}E_i \right)\\
&= -3\left( p^2+p+1+\sum_{\sigma \in F\backslash\{\vec{0},\vec{e}_1\}}(q+m_{\sigma})\right)+(p+1)(p^2+p+1)\\
&= p^3+O(p^2).
\end{align*}

\noindent (iii) Similarly, we calculate
\begin{align*}
\left( \sum_{\sigma \in G} D_{\sigma} \right)^2 &= \left( p^2+p+1+\sum_{\sigma \in F\backslash\{\vec{0},\vec{e}_1\}}(q+m_{\sigma})\right)^2-(p+1)^2(p^2+p+1)\\
&=(p^4+2p^3+O(p^2))-(p^4+3p^3+O(p^2))\\
&=-p^3+O(p^2).
\end{align*}

\noindent (iv) Lastly, by equation \eqref{eq:2} we have
\begin{align*}
\sum_{\chi \in G^{\vee}} L_{\chi}& \cdot (L_{\chi}+K_S)=\sum_{\gamma \in G} L_{\chi_{\gamma}}\cdot (L_{\chi_{\gamma}}+K_S)\\
&=\frac{1}{q^2}\sum_{\gamma \in G}[ -3q\gamma_1 (p^2+p+1)+\gamma_1^2(p^2+p+1)^2\\
&\quad +2\gamma_1(p^2+p+1)f(\gamma,\sigma)-3qf(\gamma,\sigma)+f(\gamma,\sigma)^2\\
&\quad +\gamma_1(p+1)(q-\gamma_1(p+1))(p^2+p+1)]\\
&=\frac{1}{q^2}\sum_{\gamma \in G}\left[\gamma_1^2(p^2+p+1)^2+\gamma_1(p+1)(q-\gamma_1(p+1))(p^2+p+1)\right]+O(p^2)\\
&=\frac{1}{q^2}\left[ (p^4+2p^3)\sum_{\gamma \in G}\gamma_1^2 +p^3q\sum_{\gamma \in G}\gamma_1-(p^4+3p^3)\sum_{\gamma \in G}\gamma_1^2\right]+O(p^2)\\
&=\frac{1}{q^2}\left[ -\sum_{\gamma \in G}\gamma_1^2+q\sum_{\gamma \in G}\gamma_1\right]p^3+O(p^2)\\
&=\left[-\frac{1}{q^2}\sum_{0\leq \gamma_1,\ldots,\gamma_n\leq q-1}\gamma_1^2+\frac{1}{q}\sum_{0\leq \gamma_1,\ldots,\gamma_n\leq q-1}\gamma_1\right]p^3+O(p^2)\\
&=\left[-\frac{1}{q^2}\left(q^{n-1}\frac{q(q-1)(2q-1)}{6}\right)+\frac{1}{q}\left(q^{n-1}\frac{q(q-1)}{2}\right)\right]p^3+O(p^2)\\
&=\left[\frac{1}{6}q^{n-2}(q^2-1)\right]p^3+O(p^2).
\end{align*}

\end{proof}

We're now in a position to quickly prove the proposition.
\begin{proof}
Recall that
\begin{align*}
K_{\tilde S}^2&=q^n\left( K_S+\frac{q-1}{q}\sum_{\sigma}D_{\sigma}\right)^2.\\
\intertext{Applying parts (i)-(iii) of Lemma 3 then immediately yields}
K_{\tilde S}^2&=q^{n-2}(q^2-1)p^3+O(p^2).
\intertext{Similarly, we have}
\chi ({\mathcal O}_{\tilde S})&=q^n\chi ({\mathcal O}_S)+\frac{1}{2}\sum_{\chi}L_{\chi}\cdot (L_{\chi}+K_S)\\
&=\frac{1}{12}q^{n-2}(q^2-1)p^3+O(p^2).
\end{align*}
The result follows.
\end{proof}

The following is now immediate.

\begin{corollary*}For each prime $q\geq 3$ and all sufficiently large primes $p\equiv -1 \pmod q$, there exists a surface $\tilde S$ in characteristic $p$, nonsingular and of general type, with $\frac{K_{\tilde S}^2}{\chi ({\mathcal O}_{\tilde S})}>9$.
\end{corollary*}

\paragraph{\bf{Future Considerations:}}
R. Pardini has posed several interesting questions, regarding the minimality of these surfaces and variations of the construction.  The minimality of the surfaces remains to be seen.  As for the construction, the
requirement $p\equiv -1 \pmod q$ seems quite arbitrary.  It is
likely that variations of this construction could produce similar
surfaces for arbitrary primes $p$.

\end{document}